\documentclass[a4paper,10pt]{article}
\usepackage{amsmath,amssymb}
\usepackage{lineno,hyperref}
\newtheorem{theorem}{Theorem}[section]

\newtheorem{remark}{Remark}[section]

\newtheorem{example}{Example}[section]

\begin{document}

\newcommand{\stt}{\stackrel{\longrightarrow}}
\newcommand{\st}{\stackrel{\to}}
\newcommand{\sr}{\stackrel{\rule{5pt}{0,5pt}}}
\newcommand{\srr}{\stackrel{\rule{9pt}{0,5pt}}}
\newcommand{\di}{\displaystyle}
\newcommand{\mm}{\medskip}
\newcommand{\ov}{\over}
\newcommand{\ld}{\ldots}
\newcommand{\wh}{\widehat}
\newcommand{\cd}{\cdots}
\newcommand{\te}{\theta}
\newcommand{\ty}{\infty}
\newcommand{\pp}{\prime}
\newcommand{\lft}{\Leftrightarrow}
\newcommand{\sm}{\smallskip}
\newcommand{\ri}{\Rightarrow}
\newcommand{\fo}{\forall}
\newcommand{\ga}{\gamma}
\newcommand{\Ga}{\Gamma}
\newcommand{\si}{\sigma}
\newcommand{\Si}{\Sigma}
\newcommand{\Om}{\Omega}
\newcommand{\su}{\subset}
\newcommand{\de}{\delta}
\newcommand{\ep}{\varepsilon}
\newcommand{\va}{\varphi}
\newcommand{\qu}{\quad}
\newcommand{\la}{\lambda}
\newcommand{\al}{\alpha}
\newcommand{\be}{\beta}
\newcommand{\pa}{\partial}
\newcommand{\un}{\underline}
\newcommand{\ol}{\overline}
\newcommand{\noa}{\noalign{\medskip}}
\newcommand{\tg}{\hbox{tg}\,}
\newcommand{\grad}{\hbox{grad}\,}
\newcommand{\sh}{\hbox{sh}\,}
\newcommand{\ch}{\hbox{ch}\,}
\newcommand{\sgn}{\hbox{sgn}}
\newcommand{\rar}{\rightarrow}
\newcommand{\ctg}{\hbox{ctg}\,}
\newcommand{\br}{\mathbb R}
\newcommand{\bz}{\mathbb Z}
\newcommand{\bn}{\mathbb N}
\newcommand{\bc}{\mathbb C}
\newcommand{\bq}{\mathbb Q}
\newcommand{\bm}{\mathbb M}
\newcommand{\arctg}{\hbox{arctg}\,}
\newcommand{\wt}{\widetilde}

\baselineskip 12.5pt
\title{\LARGE \bf Optimal control problems for stress tensor in plastic plane medium}

\author{\begin{minipage}[t]{4.3in} \normalsize \baselineskip 12.5pt
\centerline{SIMONA DINU and ANDREEA BEJENARU}
\vspace{1cm}
\centerline{University Politehnica of Bucharest}
\centerline{Faculty of Applied Sciences}
\centerline{Department of Mathematics-Informatics}
\centerline{Splaiul Independentei 313}
\centerline{060042 Bucharest, Romania}
\centerline{simogrigo@yahoo.com; bejenaru.andreea@yahoo.com}
\end{minipage} \kern 0in
\\ \\ \hspace*{-10pt}}
\maketitle 
\begin{abstract}
This paper joins some concepts from Mechanics, Partial Differential Equations and Control Theory in order
to solve bi-time optimization problems related to stress tensor in plastic deformations. The main goal is to analyze some optimal control problems constrained by the equilibrium equations of the stress tensor in perfect plastic plane medium. As consequence of this approach,  a natural split of the constraints arises, leading to integrability conditions and changes a classical variational problem into an optimal control one. The final outcomes confirm all the expectations related to the physical features of plastic deformations phenomenon. \\
{\bf Keywords:} multitime maximum principle, complete integrability conditions, non-linear PDE system, perfect plastic plane medium.\\
{\bf Mathematics Subject Classification (2000):} 70H06, 70H30, 70Q05, 49J20.
\end{abstract}



\section {Introduction} Over the last few years, a research team from University Politehnica of Bucharest, supervised by professor dr. C. Udriste, has proved some very interesting theoretical facts related to multitime optimal control (\cite {UMP}-\cite{Ump}, \cite{UDT}-\cite{[25]}). Lately, in order to overcome the theoretical dimension and to achieve some practical confirmation, the team focused on applying these results for meaningful problems from different scientific areas: differential geometry, statistics, mechanics etc. (\cite{[19]},\cite {[26]},\cite{[28]}). This paper is the consequence of these efforts, proving the applicability and utility of a general multitime maximum principle and also emphasizing remarkable particularities of the specific mechanical problem approached here.

Section 2 analyzes the complete integrability conditions for a quasi-linear (non-elementary) PDE system with two variables (a bi-time) and a multi-dimensional undetermined map (2-sheet state variable). The result obtained here is applied in order to describe the integrability context for the PDEs defining the 2-dimensional stress tensor. The major outcome of the section consists in emphasizing a natural split of a non-linear PDE system, over state gradients and a natural insertion of control variables, leading to manageable integrability conditions.
Section 3 describes optimization processes constrained by non-linear PDEs,  a bi-time maximum principle. The main result is adapted for constrained variational problems. Moreover, this section points out the natural transformation of variational processes to optimal control ones, via the canonical controls resulted from integrability requirements.   Finally, Section 4 applies the theoretical facts from the first two sections in order to describe the solution of a variational process involving the stress tensor in perfect plastic plane medium.

\section{Complete integrability conditions\\
 for plane quasi-linear PDE systems}

The main goal of this Section is to phrase the complete integrability conditions for the PDE system describing the stress tensor in perfect plastic medium. Basically, this starts with finding a pertinent description for integrability conditions for an abstract, general problem. Three important consequences shall derive from our attempt: 1) there is a natural split of an arbitrary quasi-linear PDE system, separating state gradients one from each other, generating manageable integrability conditions, 2) this natural split is determined by a natural insertion of control variables and 3) any variational problem constrained by an arbitrary linear PDE system may be naturally rephrased as an optimal control problem via the control variables mentioned above.

\subsection{ One state-variable plane PDE systems}

Let $\Omega$ be a compact subset of $ R^2$, with global coordinates $t=(t^1,t^2)$ and let $\Sigma=\pa\Omega$ denote its boundary. If $A=(A^{\alpha\beta}:\Omega\times R\times R^k\rightarrow R)_{\al,\be=\overline{1,2}}$ and $B=(B^{\alpha}:\Omega\times R\times R^k\rightarrow R)_{\al=\overline{1,2}}$ denote some $2\times 2$, respectively $2\times 1 $ tensors,  they define the following \textit{one state quasi-linear PDE system} (\cite{CR}):
$$A(t,x(t),u(t))\grad x(t)=B(t,x(t),u(t)),\leqno(1)$$
where, for the sake of simplicity, $t$ stands for $(t^1,t^2)$ and denotes the {\it bi-time variable}, $x$ represents the {\it state variable} (unknown differentiable two-sheet, involved in the PDEs via its gradient) and $u(t)=(u^k(t))_{k=\overline{1,N}}$ denotes the {\it control variable} (parameter for the PDE system). Written explicitly, the PDEs are:
$$
\left\{\begin{array}{ll}
A^{11}(t,x(t),u(t))\di\frac{\pa x}{\pa t^1}+A^{12}(t,x(t),u(t))\di\frac{\pa x}{\pa t^2}=B^1(t,x(t),u(t));\\ \noa
A^{21}(t,x(t),u(t))\di\frac{\pa x}{\pa t^1}+A^{22}(t,x(t),u(t))\di\frac{\pa x}{\pa t^2}=B^2(t,x(t),u(t)).
\end{array}\right.\leqno(1')$$

Suppose $x=x(t^1,t^2)$ is a solution for the foregoing PDE system, defined implicitly by $\Phi(t,x(t))=0$.

Applying the Implicit Functions Theorem, it follows
$$
\di\frac{\pa x}{\pa t^\alpha}=-\di\frac{\di\frac{\pa \Phi}{\pa t^\alpha}}{\di\frac{\pa \Phi}{\pa x}},\,\,\alpha=\overline{1,2}.
$$
Substituting these in $(1')$ leads to a linear PDE

$$
\left\{\begin{array}{ll}
A^{11}\di\frac{\pa \Phi}{\pa t^1}+A^{12}\di\frac{\pa \Phi}{\pa t^2}+B^1\di\frac{\pa \Phi}{\pa x}=0\\ \noa
A^{21}\di\frac{\pa \Phi}{\pa t^1}+A^{22}\di\frac{\pa \Phi}{\pa t^2}+B^2\di\frac{\pa \Phi}{\pa x}=0
\end{array}\right.,
$$
meaning $\grad \Phi \bot \overline{A}^{\alpha},\,\,\alpha=\overline{1,2},$ where $\overline{A}^{\alpha}=(A^{\alpha 1},A^{\alpha 2},B^{\alpha})$. Since $\grad\Phi$ and $\overline{A^1}\times \overline{A^2}$ are collinear,  {\it the vectorial expression for complete integrability conditions for plane non-linear PDE systems} is derived:
$$\grad\Phi\times(\overline{A^1}\times \overline{A^2})=0.\leqno(2)$$
Explicitly, this means
$$
\left\{\begin{array}{cc}
R\di\frac{\pa x}{\pa t^2}+Q&=0\\ \noa
R\di\frac{\pa x}{\pa t^1}+P&=0\\ \noa
-Q\di\frac{\pa x}{\pa t^1}+P\di\frac{\pa x}{\pa t^2}&=0,
\end{array}\right.
$$
where $(P,Q,R)=\overline{A^1}\times \overline{A^2}=(A^{12}B^{2}-A^{22}B^1,A^{21}B^1-A^{11}B^2,A^{11}A^{22}-A^{12}A^{21})$. The third relation above is a consequence of the previous ones and the differentiation of the first and the second one with respect to $t^1$, respectively $t^2$  is leading to {\it the explicit integrability condition for one state plane non-linear PDE systems}:
$$
\di\frac{\pa}{\pa t^2}\left(P-x\di\frac{\pa R}{\pa t^1}\right)=\di\frac{\pa}{\pa t^1}\left(Q-x\di\frac{\pa R}{\pa t^2}\right). \leqno(3)
$$

\subsection{ Multi-state variables plane PDE systems}

In the following, the main interest is to extend the previous result for linear PDEs having $n$ state variables. Let $A_i=(A^{\alpha\beta}_{i}:\Omega\times R^n\times R^k\rightarrow R)_{\al,\be=\overline{1,2}}$ and $B=(B^{\alpha}:\Omega\times R^n\times R^k\rightarrow R)_{\al=\overline{1,2}}$ denote some $2\times 2$, respectively $2\times 1 $ matrix fields, where $i=\overline{1,n}$ denotes the state variables index. Let us start with the following \textit{quasi-linear PDE system}:
$$A_i(t,x(t),u(t))\grad x^i(t)=B(t,x(t),u(t)).\leqno(4)$$

The basic and also an original idea  used here consists in separating the gradients of the two state variables via some additional control. More precisely, it consists in rewriting the above PDEs as {\it a split non-linear PDE system}:
$$\left\{\begin{array}{ll}A_i(t,x(t),u(t))\grad x^i(t)=v_i(t),\,\,i=\overline{1,n-1}\,\, \mbox{(no sum)};\\
A_n(t,x(t),u(t))\grad x^n(t)=B(t,x(t),u(t))-\sum_{i=1}^{n-1}v_i(t).\end{array}\right.\leqno(5)$$

Then, by applying the result obtained in the previous Section, {\it the complete integrability conditions for $n$-state variables quasi-linear PDE systems} may be phrased:
$$
\begin{array}{ll}
\di\frac{\pa}{\pa t^2}\left(P_i-x^i\di\frac{\pa R_i}{\pa t^1}\right)=\di\frac{\pa}{\pa t^1}\left(Q_i-x^i\di\frac{\pa R_i}{\pa t^2}\right),\,\,i=\overline{1,n},\\
\end{array} \leqno(6)
$$
where
$$
(P_i,Q_i,R_i)=(A^{12}_iv^{2}_i-A^{22}_iv^1_i,A^{21}_iv^1_i-A^{11}_iv^2_i,\det A_i),\,\,\forall i=\overline{1,n-1}
$$
and
$$\left\{\begin{array}{ll} P_n=A^{12}_n\left(B^{2}-\sum_{i=1}^{n-1}v^2_i(t)\right)-A^{22}_n\left(B^1-\sum_{i=1}^{n-1}v^1_i(t)\right);\\
Q_n=A^{21}_2\left(B^1-\sum_{i=1}^{n-1}v^1_i(t)\right)-A^{11}_2\left(B^{2}-\sum_{i=1}^{n-1}v^2_i(t)\right);\\
R_n=\det A_n.\end{array}\right.
$$

\subsection{The PDE system of stress tensor in
 perfect plastic medium}

As  mentioned at the very beginning of this section, the content of the paper is related to the stress tensor for deformations in perfect plastic plane medium. In this paragraph,  a PDE system associated to this geometric object is defined and rewritten in a more manageable way and, by applying the theoretical results obtained above, the corresponding complete integrability conditions are phrased. For this, it makes total use of the inspired idea of splitting linear PDE systems over gradients; this separation will also prove to be of major utility in the further development of the paper.

There is no novelty that the stress tensor components for perfect plastic plane medium are described by a constrained non-linear PDE system (equilibrium condition, see \cite{Bol}, \cite{NT}):
$$
\left\{\begin{array}{l}
\di\frac{\pa\sigma_{xx}}{\pa x}+\di\frac{\pa \sigma _{xy}}{\pa y}=0,\\ \noa
\di\frac{\pa\sigma _{xy}}{\pa x}+\di\frac{\pa\sigma _{yy}}{\pa y}=0,\\ \noa
(\sigma _{yy}-\sigma _{xx})^2+4\sigma _{xy}^2=4K^2,
\end{array}\right.\leqno(7)
$$
where $\sigma _{xx}$, $\sigma _{yy}$ (normal stresses) and $\sigma _{xy}$ (shear stress) denote the components of the symmetric stress tensor and $K$ is the radius of the Mohr's circle (depending, as Mohr's equation proves it, on the average of the normal stresses).

Using a classic and natural change for state variables (meaning a rotation of angle $\frac{\va}{2}$ in order to replace the given coordinate system with principal directions) it leads to $\sigma _{xx}=\rho -K\cos\va,\,\,\sigma _{yy}=\rho +K\cos\va,\,\,\sigma _{xy}=K\sin\va$ and
may be written {\it the polar non-linear PDE system}:
$$
\left\{\begin{array}{l}
\di\frac{\pa \rho }{\pa x}-\frac{\pa K}{\pa x}\cos\va+K\sin\va\di\frac{\pa \va}{\pa x}+\frac{\pa K}{\pa y}\sin\va+K\cos\va\di\frac{\pa \va}{\pa y}=0,\\ \noa
\di\frac{\pa \rho}{\pa y}+\frac{\pa K}{\pa y}\cos\va-K\sin \va\di\frac{\pa \va}{\pa y}+\frac{\pa K}{\pa x}\sin\va+K\cos\va\di\frac{\pa \va}{\pa x}=0.
\end{array}\right.\leqno ( 8)
$$

It follows that $(x,y)$ may be identified with the bi-time variable, $(\rho,\va,K)$ is the three-dimensional state variable and the matrices of the non-linear PDE system are
$$
A_1=\left(\begin{array}{cc}
          1 & 0  \\
         0& 1 \end{array}\right),\,\,
A_2=\left(\begin{array}{cc}
 -\cos\va & \sin\va \\
 \sin\va & \cos\va \end{array}\right),\,\,
 A_3=\left(\begin{array}{cc}
 K\sin\va & K\cos\va \\
 K\cos\va &-K\sin\va \end{array}\right).
$$

 It is obvious that initial control variables are missing, but some natural control variables emerge as consequence of the canonical split of the system over gradients. Applying this technique here leads to

$$
\left\{\begin{array}{l}
\di\frac{\pa \rho }{\pa x}=u,\qu \di\frac{\pa \rho}{\pa y}=v,\\ \noa
-\di\frac{\pa K }{\pa x}\cos \va+\di\frac{\pa K }{\pa y}\sin \va=\mu,\\
 \di\frac{\pa K }{\pa x}\sin \va+ \di\frac{\pa K}{\pa y}\cos\va=\nu,\\ \noa
K\sin\varphi\di\frac{\pa \varphi}{\pa x}+K\cos\varphi\di\frac{\pa \varphi}{\pa y}=-u-\mu,\\ \noa
K\cos\varphi\di\frac{\pa \varphi}{\pa x}-K\sin \varphi\di\frac{\pa \varphi}{\pa y}=-v-\nu,
\end{array}\right.\leqno (9)
$$
or, equivalent
$$
\left\{\begin{array}{l}
\di\frac{\pa \rho }{\pa x}=u,\qu \di\frac{\pa \rho}{\pa y}=v,\\ \noa
\di\frac{\pa K }{\pa x}=-\mu \cos \va +\nu \sin\va, \qu \di\frac{\pa K }{\pa y}=\mu \sin \va+\nu \cos\va, \\ \noa
K\di\frac{\pa \varphi}{\pa x}=-(u+\mu)\sin\va-(v+\nu)\cos \va,\\
 K\di\frac{\pa \varphi}{\pa y}=-(u+\mu)\cos\va+(v+\nu)\sin \va.
\end{array}\right.
$$
Finally, the {\it integrability conditions for perfect plastic medium problem} are the result of relations $(6)$ and have the simplest expression

$$
\left\{\begin{array}{ll}
\di\frac{\pa u}{\pa y}=\di\frac{\pa v}{\pa x};\\ \noa
\di\frac{\pa}{\pa y}(-\mu \cos\varphi+\nu \sin\varphi)=\di\frac{\pa}{\pa x}(\mu \sin\varphi+ \nu\cos\varphi);\\ \noa
\di\frac{\pa}{\pa y}\Big((u+\mu)\sin\va+(v+\nu)\cos \va\Big)+\frac{\pa K}{\pa y}\frac{\pa\va}{\pa x}=\di\frac{\pa}{\pa x}\Big((u+\mu)\cos\va-(v+\nu)\sin \va\Big)\\
+\di\frac{\pa K}{\pa x}\frac{\pa\va}{\pa y}.
 \end{array}\right.\leqno( 10)
$$

\begin{remark} {\it As direct consequence for the above integrability conditions, for each predefined state $\va=\va(x,y)$, the rest of the state variables $\rho$ and $K$ are solutions for the complete integrable PDE system
$$
\left\{\begin{array}{l}
\di\frac{\pa \rho }{\pa x}=\frac{\pa }{\pa x}(K\cos\va)-\frac{\pa }{\pa y}(K\sin\va);\\ \noa
\di\frac{\pa \rho}{\pa y}=-\frac{\pa }{\pa x}(K\sin\va)-\frac{\pa}{\pa y}(K\cos\va);\\ \noa
\di2\frac{\pa^2}{\pa x\pa y}(K\cos \va)=\di\frac{\pa^2}{\pa y^2}(K\sin \va)-\di\frac{\pa^2}{\pa x^2}(K\sin \va).
\end{array}\right.
$$
} \end{remark}

\section{Optimal control problems constrained\\ by
 non-linear PDE systems}

The aim of this section is to use arbitrary non-linear PDE systems as constraints for optimizing cost functionals defined as multiple integrals. In the most general framework, the main interest consists in finding
$$
\di\max_{u(\cdot)}I(u(\cdot))=\di\int_{\Omega}X(t,x(t),u(t))dv+\di\int_{\Sigma=\pa\Omega}g(t,x(t))d\sigma,\leqno(11)
$$
constrained by
$$
A_i^{\be\al}(t,x(t),u(t))\di\frac{\pa x^i}{\pa t^\al}(t)=B^{\be}(t,x(t),u(t)),\,\,\forall \be=\overline{1,2}.
$$

In the spirit of the breakthrough in the first Section, it is of major utility to replace the previous constraints with a split $n$-state variable quasi-linear PDE system
$$\left\{\begin{array}{ll}A_i(t,x(t),u(t))\grad x^i(t)=v_i(t),\,\,i=\overline{1,n-1}\,\, \mbox{(no sum)};\\
A_n(t,x(t),u(t))\grad x^n(t)=B(t,x(t),u(t))-\sum_{i=1}^{n-1}v_i(t).\end{array}\right.\leqno(12)$$

The complete integrability conditions $(6)$ define the {\it set of admissible controls} (with an initial component $u=(u^1,...,u^N)$ and some additional ones $v_i=(v_i^1,v_i^2),\,\,i=\overline{1,n-1}$):
$${\mathcal U}=\Big\{\overline{u}=(u,v):\Omega\rightarrow R^{N+2n-2}\Big |\, \overline{u} \mbox{ constrained by } (CIC_2)\Big\}$$

\begin{remark} {\it Any variational problem
$$
\di\max_{x(\cdot)}I(x(\cdot))=\di\int_{\Omega}X(t,x(t))dv+\di\int_{\Sigma=\pa\Omega}g(t,x(t))d\sigma,
$$
constrained by
$$
A_i^{\be\al}(t,x(t))\di\frac{\pa x^i}{\pa t^\al}(t)=B^{\be}(t,x(t)),\,\,\al,\be=\overline{1,2},\,\,i=\overline{1,n}
$$
becomes, after the natural split of the PDE system and the natural addition of canonical control variables, an optimal control problem (with no initial controls). This fact was anticipated in many papers describing extensions from calculus of variations to optimal control (see  \cite{Ev}, \cite{PW}, \cite{UMP}, \cite{UN}, \cite{UDT}, \cite{UM},\cite{[24]}). Nevertheless, the originality of this paper consists in the fact that explains, for the first time, the natural process of inserting canonical control variables compatible with the constraints.}
\end{remark}
In order to solve the optimal control problem and following the classical pattern, some  $C^1$ Lagrange multipliers $p^i(t)=(p^i_\al(t))_{\al=\ol{1,2}},\,i=\ol{1,n}$ are introduced in order to define the Hamiltonian
\begin{eqnarray*}
H(t,x(t),u(t),p(t))&=&X(t,x(t),u(t))+\sum_{i=1}^{n-1}p^i_\beta(t)v^{\be}_i(t)\\
&+&p^n_\beta(t)\left(B^{\be}(t,x(t),u(t))-\sum_{i=1}^{n-1}v^{\be}_i(t)\right).
\end{eqnarray*}

The initial variational problem is changed into an optimal control one, namely
$$
\di\max_{u(\cdot)}I(u(\cdot))=\di\int_{\Omega}\mathcal{L}(t,x(t),\overline{u}(t),p(t))dv+
\di\int_{\Sigma}g(t,x(t))d\sigma,\leqno(13)
$$
or
\[\max_{u(\cdot)}I(u(\cdot))=\int_{\Omega}\left[H-p^i_\beta A_i^{\be\al}\frac{\pa x^i}{\pa t^\al}\right]dv+\int_{\Sigma}g(t,x)d\sigma.\]

Although several more general approaches on control theory have been analyzed so far (see \cite{Zu}, \cite{Ubang}), this paper bases on the assumption that the optimal problem $((12),(13))$ admits a continuous optimal control $\overline{u}^*(t)=(u^*(t),v^*(t))\in Int({\mathcal U})$ which generates an optimal state $x^*(t)$. Let us consider a variation of the control $\overline{u}_\xi(t)=\overline{u}^*(t)+\xi h(t)$, where $h$ is an arbitrary continuous vector. Since $\overline{u}^*(t)\in Int({\mathcal U})$  and a continuous function  over a compact set $\Omega$  is bounded, there exists $\xi_h>0$ such that
$\overline{u}_{\xi}(t)=\in Int\mathcal{U}(t),\,\forall |\xi|<\xi_h.$ This $\xi$ is used in our arguments. Furthermore, let $x^i_{\xi}(t)$ be the corresponding variation of the optimal state and $y(t)=\di\frac{\pa x_\xi}{\pa \xi}(t)\Big|_{\xi=0}$ be the variation vector field.
For $|\xi|<\xi_h$, let us define the function:
$$
\begin{array}{ll}
I(\xi)&=\di\int_{\Omega}\mathcal{L}(t,x^i_{\xi}(t),\overline{u}_{\xi}(t),p(t))dv+
\di\int_{\Sigma}g(t,x^i_{\xi}(t))d\sigma=\\ \noa
&=\int_{\Omega}\left[H(t,x^i_{\xi}(t),\overline{u}_{\xi}(t),p(t))-p^i_\beta(t) A_i^{\be\al}(t,x^i_{\xi}(t),\overline{u}_{\xi}(t))\frac{\pa x^i_\xi}{\pa t^\al}(t)\right]dv+\\ \noa
&+\int_{\Sigma}g(t,x_{\xi}(t))d\sigma.
\end{array}
$$

Differentiating the function $I(\xi)$ leads to
$$
\begin{array}{ll}
&I'(\xi)=\di\int_{\Omega}\left\{\di\frac{\pa H}{\pa x^i}(t,x_{\xi},\overline{u}_{\xi},p)\di\frac{\pa x^i_\xi}{\pa \xi}+\di\frac{\pa H}{\pa \overline{u}^a}(t,x_{\xi},\overline{u}_{\xi},p)\di\frac{\pa \overline{u}^a_\xi}{\pa \xi}\right.\\ \noa
&-p_\be^i A_i^{\be\al}(t,x_{\xi},\overline{u}_{\xi})\di\frac{\pa^2 x^i_\xi}{\pa \xi\pa t^\al}
-p_\be^i\left[\di\frac{\pa A_i^{\be\al}}{\pa x^j}(t,x_{\xi},\overline{u}_{\xi})\di\frac{\pa x^j_\xi}{\pa \xi}\right.\\ \noa
&\left.\left.+\di\frac{\pa A_i^{\be\al}}{\pa \overline{u}^a}(t,x_{\xi},\overline{u}_{\xi})\di\frac{\pa \overline{u}^a_\xi}{\pa \xi}\right]\di\frac{\pa x^i_{\xi}}{\pa t^\al}\right\}dv+\di\int_{\Sigma}\frac{\pa g}{\pa x^i}(t,x_\xi)\frac{\pa x^i_\xi}{\pa \xi}d\sigma.
\end{array}
$$

For $\xi=0$, it follows
$$
\begin{array}{ll}
&I'(0)=\di\int_{\Omega}\left\{\di\frac{\pa H}{\pa x^i}(t,x^*,\overline{u}^*,p)y^i+\di\frac{\pa H}{\pa \overline{u}^a}(t,x^*,\overline{u}^*,p)h^a-p_\be^iA_i^{\be\al}(t,x^*,\overline{u}^*)\di\frac{\pa y^i}{\pa t^\al}\right.\\ \noa
&\left.-p_\be^i\Big[\di\frac{\pa A_i^{\be\al}}{\pa x^j}(t,x^*,\overline{u}^*)y^j+\di\frac{\pa A_i^{\be\al}}{\pa \overline{u}^a}(t,x^*,\overline{u}^*)h^a\Big]\di\frac{\pa x^{i*}}{\pa t^\al}\right\}dv+\di\int_{\Sigma}\frac{\pa g}{\pa x^i}(t,x_\xi)y^i d\sigma.
\end{array}
$$

By substituting the term
$-p_\be^iA_i^{\be\al}\di\frac{\pa y^i}{\pa t^\al}$ with
$-\di\frac{\pa }{\pa t^\al}\left[p_\be^iA_i^{\be\al}y^i\right]+\di\frac{\pa }{\pa t^\al}(p_\be^iA_i^{\be\al})y^i,$ the relation above may be rewritten
$$
\begin{array}{ll}
&I'(0)=+\di\int_{\Sigma}\left[\frac{\pa g}{\pa x^i}y^i-p_\be^i A_i^{\be\al}y^i n^\al\right] d\sigma\\
&+\di\int_{\Omega}\left\{\di\frac{\pa H}{\pa x^i}y^i+\di\frac{\pa H}{\pa \overline{u}^a}h^a+\di\frac{\pa }{\pa t^\al}\Big[p_\be^iA_i^{\be\al}\Big]y^i-p_\be^i\Big(\di\frac{\pa A_i^{\be\al}}{\pa x^j}y^j+\di\frac{\pa A_i^{\be\al}}{\pa \overline{u}^a}h^a\Big)\di\frac{\pa x^{i*} }{\pa t^\al}\right\}dv\\ \noa

&= \di\int_{\Omega}\left[\di\frac{\pa H}{\pa x^i}+\di\frac{\pa }{\pa t^\al}\Big(p_\be^iA_i^{\be\al}\Big)-p_\be^j\di\frac{\pa A_j^{\be\al}}{\pa x^i}\di\frac{\pa x^{j*} }{\pa t^\al}\right]y^i+\left[\di\frac{\pa H}{\pa \overline{u}^a}-p_\be^i\di\frac{\pa A_i^{\be\al}}{\pa \overline{u}^a}\di\frac{\pa x^{i*} }{\pa t^\al}\right]h^a dv\\ \noa

&+\di\int_{\Sigma}\left[\frac{\pa g}{\pa x^i}-p_\be^i A_i^{\be\al} n_\al\right]y^i d\sigma.\\

\end{array}
$$

Since $\xi=0$ is the critical point of $I(\xi)$, it follows $I'(0)=0$. Defining the set ${\mathcal P}=\{p^*(t)\}$ of optimal costate variables as the set of solutions for the Cauchy problem
$$
\left\{\begin{array}{l}
\di\frac{\pa H}{\pa x^i}(t,x^*,\overline{u}^*,p^*)+\di\frac{\pa}{\pa t^\al}\Big(p_\be^{i*}A_i^{\be\al}(t,x^*,\overline{u}^*)\Big)-p^{j*}_\be(t)\di\frac{\pa A_j^{r\al}}{\pa x^i}(t,x^*,\overline{u}^*)\di\frac{\pa x^{j*}}{\pa t^\al}(t)=0,\\ \noa
\di\frac{\pa g}{\pa x^i}(t,x^*)-p_\be^{i*}(t)A_i^{\be\al}(t,x^*,\overline{u}^*)n_\al(t)\Big|_{\Sigma}=0,\,\,\forall i=\overline{1,n} \,\,(\mbox{ i.e. no sum over } i),
\end{array}\right.
$$
leads to the optimum critical point condition
$$
\di\frac{\pa H}{\pa \overline{u}^a}(t,x^*,\overline{u}^*,p^*)-p_\be^{i*}(t)\di\frac{\pa A_i^{\be\al}}{\pa \overline{u}^a}(t,x^*,\overline{u}^*)\di\frac{\pa x^{*i}}{\pa t^\al}=0
$$

According to the above statements, the following result may be phrased:
\begin{theorem}
If $\overline{u}^*(t)$ is an optimal solution for the problem $((12),(13))$ and $x^*(t)$ is the corresponding optimal state, then there exist the costate vectors $p^i=(p_\be^i)_{\be=\ol{1,2}}$ such that $(x^*,\overline{u}^*,p^*)$ satisfies
$$
\left\{\begin{array}{ll}
v_i^{\be *}=\di A_i^{\be\al}(t,x^*,\ol{u}^*)\di\frac{\pa x^{i*}}{\pa t^\al}\,\,\mbox{(no sum over i)};\\ \noa
B^\be(t,x^*,\ol{u}^*)-\sum_{i=1}^{n-1}v_i^{\be *}=\di A_n^{\be\al}(t,x^*,\ol{u}^*)\frac{\pa x^{n*}}{\pa t^\al},\end{array}\right.\leqno(14)
$$

$$
\begin{array}{ll}
\di\frac{\pa H}{\pa x^i}(t,x^*,\ol{u}^*,p^*)+\di\frac{\pa}{\pa t^\al}\Big(p_\be^{i*}A_i^{\be\al}(t,x^*,\ol{u}^*)\Big)-p^{j*}_\be(t)\di\frac{\pa A_j^{\be\al}}{\pa x^i}(t,x^*,\overline{u}^*)\di\frac{\pa x^{j*}}{\pa t^\al}=0\\
\noa\mbox{(no sum over i),}
\end{array}\leqno(15)
$$
$$
\di\frac{\pa}{\pa t^2}\left(P_i-x^i\di\frac{\pa R_i}{\pa t^1}\right)=\di\frac{\pa}{\pa t^1}\left(Q_i-x^i\di\frac{\pa R_i}{\pa t^2}\right);\\
 \leqno(16)
$$

$$
\di\frac{\pa H}{\pa \ol{u}^a}(t,x^*,\ol{u}^*,p^*)-p_\be^{i*}(t)\di\frac{\pa A_i^{\be\al}}{\pa \ol{u}^a}(t,x^*,\ol{u}^*)\di\frac{\pa x^{*i}}{\pa t^\al}(t)=0,\,\,\forall a=\overline{1,N+2n-2};\leqno(17)
$$

$$
\di\frac{\pa g}{\pa x^i}(t,x^*)-p_\be^{i*}A_i^{\be\al}(t,x^*,\ol{u}^*)n_\al\Big|_{\Sigma}=0 \,\,\mbox{(no sum over i)}.
\leqno(18)
$$
\end{theorem}

\begin{remark}{\it For a variational problem as the one described in Remark 2.1, the relations above gain a much simple expression and they are much more easy to work with, due to the absence of initial control variables. More precisely, $\ol{u}$ consists only in the canonic additional control variable $v$. Then,
$$
\left\{\begin{array}{ll}
v_i^{\be *}=\di A_i^{\be\al}(t,x^*)\frac{\pa x^{*i}}{\pa t^\al} \,\,\mbox{(no sum over i)};\\
B^\be(t,x^*)-\sum_{i=1}^{n-1}v_i^{\be *}=\di A_n^{\be\al}(t,x^*)\frac{\pa x^{*n}}{\pa t^\al},
\end{array}\right.\leqno(19)
$$

$$\begin{array}{c}
\di\frac{\pa H}{\pa x^i}(t,x^*,v^*,p^*)+\di\frac{\pa p_\be^{i*}}{\pa t^\al}A_i^{\be\al}(t,x^*)+p_\be^{i*}\di\frac{\pa A_i^{\be\al}}{\pa t^\al}(t,x^*)\\
-\di\sum_{j\neq i}p^{j*}_\be(t)\di\frac{\pa A_j^{\be\al}}{\pa x^i}(t,x^*)\di\frac{\pa x^{*j}}{\pa t^\al}=0$$ \,\,\mbox{(no sum over i)};\end{array}
\leqno(20)
$$
$$
\di\frac{\pa H}{\pa v^\al}(t,x^*,v^*,p^*)=0,\,\,\forall \al=\ol{1,2};\leqno(21)
$$

$$
\di\frac{\pa g}{\pa x^i}(t,x^*)-p_\be^{i*}A_i^{\be\al}(t,x^*)n_\al\Big|_{\Sigma}=0\,\,\mbox{(no sum over i)}.
\leqno(22)
$$
}
\end{remark}
\section{An optimal control problem constrained\\ by
 perfect plane medium PDEs}

This section combines the results related to integrability of the PDEs describing the stress tensor components in perfect plane medium (Section 2) with the maximum principle derived in Section 3, in order to completely determine a solution for a given optimization problem.
Similar ideas have successfully been applied  before in connection with stochastic theory, differential geometry, deformation theory, electromagnetic fields (\cite{Gia},  \cite{[19]},  \cite {[26]}, \cite{UTR}, \cite{[25]}, \cite{[28]}). It is emphasized once more, if necessary, the high applicability of the optimal control techniques in mechanical processes.

\begin{example}{\rm Let $\Omega=D_2$ denote the unit disc in $R^2$ and let $\Sigma=S_2$ be its boundary (the unit circle). Let us consider the following variational problem:
$$
\max_{(\va(\cdot))}\di\int_{S_2}\va(x,y)d\si,\leqno(23)
$$
constrained by PDE system:
$$
\left\{\begin{array}{l}
\di\frac{\pa \rho }{\pa x}-\frac{\pa K}{\pa x}\cos\va+K\sin\va\di\frac{\pa \va}{\pa x}+\frac{\pa K}{\pa y}\sin\va+K\cos\va\di\frac{\pa \va}{\pa y}=0,\\ \noa
\di\frac{\pa \rho}{\pa y}+\frac{\pa K}{\pa y}\cos\va-K\sin \va\di\frac{\pa \va}{\pa y}+\frac{\pa K}{\pa x}\sin\va+K\cos\va\di\frac{\pa \va}{\pa x}=0.
\end{array}\right.\leqno(24)
$$

The natural approach of this problem stands on the separation of the state gradients, via some additional, canonic control variables. Denoting these canonic control variables with $u(x,y)$, $v(x,y)$, $\mu(x,y)$ and $\nu(x,y)$, the variational problem
may be reconsidered as an optimal control problem constrained by:
$$
\left\{\begin{array}{l}
\di\frac{\pa \rho }{\pa x}=u,\qu \di\frac{\pa \rho}{\pa y}=v,\\ \noa
\di\frac{\pa K }{\pa x}=-\mu \cos \va +\nu \sin\va, \qu \di\frac{\pa K }{\pa y}=\mu \sin \va+\nu \cos\va, \\ \noa
K\di\frac{\pa \varphi}{\pa x}=-(u+\mu)\sin\va-(v+\nu)\cos \va,\\
 K\di\frac{\pa \varphi}{\pa y}=-(u+\mu)\cos\va+(v+\nu)\sin \va.
\end{array}\right.
$$

Multiplying each relation above with a corresponding Lagrange multiplier $p_1(x,y)$, $p_2(x,y)$, $r_1(x,y)$, $r_2(x,y)$, $q_1(x,y)$, $q_2(x,y)$ respectively, the corresponding Hamiltonian writes
$$\begin{array}{c}
H=p_1u+p_2v
+r_1\Big(-\mu\cos\va +\nu\sin\va \Big)+r_2\Big(\mu\sin\va +\nu\cos\va \Big)\\
+q_1\Big(-(u+\mu)\sin\va -(v+\nu)\cos\va \Big)+q_2\Big(-(u+\mu)\cos\va +(v+\nu)\sin\va \Big)
\end{array}.
$$

Application of the maximum principle for variational problems (Remark 2.2) generates the following PDE systems:
$$
\left\{\begin{array}{rl}
\di\frac{\pa p_1}{\pa x}+\di\frac{\pa p_2}{\pa y}&=0;\\ \noa
\di\frac{\pa r_1}{\pa x}+\di\frac{\pa r_2}{\pa y}&=\di q_1\frac{\pa \va}{\pa x}+q_2\frac{\pa \va}{\pa y}\\
\di\frac{\pa (Kq_1)}{\pa x}+\di\frac{\pa (Kq_2)}{\pa y}&=\sin\va\Big(-r_1\mu+r_2\nu -q_1(v+\nu)-q_2(u+\nu)\Big)\\ \noa
&-\cos\va\Big(r_1\nu+r_2\mu-q_1(u+\mu)+q_2(v+\nu)\Big),
\end{array}\right.\leqno(25)
$$
$$
\left\{\begin{array}{l}
 p_1=q_1\sin\va+q_2\cos\va;\,\, p_2=q_1\cos\va-q_2\sin\va; \\ \noa
 r_1=-q_2;\,\,r_2=q_1.
\end{array}\right.\leqno(26)
$$
and
$$
\left\{\begin{array}{l}
p_1n_1+p_2n_2=\di\frac{\pa g}{\pa\rho}\Rightarrow p_1x+p_2y=0,\\ \noa
q_1n_1+q_2n_2=\di\frac{\pa g}{\pa \va}\Rightarrow q_1x+q_2y=1,\\ \noa
r_1n_1+r_2n_2=\di\frac{\pa g}{\pa K}\Rightarrow r_1x+r_2y=0.
\end{array}\right. \mbox{ on } \Sigma.\leqno(27)
$$
Using relations $(26)$, system $(25)$ becomes
$$
\left\{\begin{array}{rl}
\di\frac{\pa p_1}{\pa x}+\di\frac{\pa p_2}{\pa y}&=0;\\ \noa
\di-\frac{\pa q_2}{\pa x}+\di\frac{\pa q_1}{\pa y}&=\di q_1\frac{\pa \va}{\pa x}+q_2\frac{\pa \va}{\pa y}\\
\di\frac{\pa q_1}{\pa x}+\di\frac{\pa q_2}{\pa y}&=\di-q_1\frac{\pa \va}{\pa y}+q_2\frac{\pa \va}{\pa x},
\end{array}\right.\leqno(28)
$$
with the particular admissible co-states $p_1(x,y)=y$, $p_2(x,y)=-x$, $q_1(x,y)=r_2(x,y)=y\sin\va(x,y)-x\cos\va(x,y)$ and $q_2(x,y)=-r_1(x,y)=y\cos\va(x,y)+x\sin\va(x,y)$. Introducing them into the boundary constrains (27), leads to
$$(y^2-x^2)\cos\va+2xy\sin\va=1 \mbox{ on } S_2.$$

Since $(y^2-x^2)^2+(2xy)^2=(x^2+y^2)^2=1$ to $S_2$, it follows
$$
\cos\va(x,y)=y^2-x^2,\,\, \sin\va(x,y)=2xy \mbox{ on } S_2.
$$
Extending this solution on $D_2$, gives:
$$
\cos\va(x,y)=\di\frac{-x^2+y^2}{x^2+y^2},\,\sin\va(x,y)=\di\frac{2xy}{x^2+y^2},
$$
hence $\va(x,y)=\arctg\di\frac{2xy}{y^2-x^2}$.

By introducing these in Remark 1.1 (combining integrability conditions with the initial PDE constraints, for a given state variable $\va(x,y)$), it follows that the state variable $K$ needs to satisfy the second order $PDE$
$$2(y^2-x^2)\frac{\pa^2 K}{\pa x\pa y}-2xy\Big(\frac{\pa^2 K}{\pa y^2}-\frac{\pa^2 K}{\pa x^2}\Big)+4\Big(y\frac{\pa K}{\pa x}-x\frac{\pa K}{\pa y}\Big)=0.$$

Look for particular solutions of type $K(x,y)=A(x)+B(y)$ leads to
$$K(x,y)=\al(x^2+y^2)+\be\frac{1}{x}+\ga\frac{1}{x}+\delta.$$

In particular,

\begin{enumerate}
\item if $K(x,y)=\al(x^2+y^2)$, then $\rho(x,y)=-2\al(x^2+y)$, therefore $K=-\frac{1}{2}\rho;$
\item if $K(x,y)=\be\frac{1}{x}$, then $\rho(x,y)=\be\frac{1}{x}$ and $K=\rho;$
\item similarly, for $K(x,y)=\ga\frac{1}{x}$, $\rho(x,y)=\ga\frac{1}{y}$ and $K=\rho.$
\item Finally, if $K(x,y)=\delta,$  then $\rho(x,y)=-\delta \ln(x^2+y^2).$
\end{enumerate}}

\end{example}

\section{Conclusions and later development}

 The example analyzed above explains the utility of a multitime Pontryaguin maximum principle when dealing with mechanical phenomenons described by non-linear PDEs. The spectacular outcome is the fact  that the solutions emphasize some reasonable expectation related to the dependence of $K$ on $\rho$ (fact already known as basic feature in plastic deformations).

 Moreover, the approach described in this paper encourages us to apply similar techniques for even more sophisticated problems in physics (generated, for example by  higher order PDEs) or even in differential geometry or other scientific areas.

\section{ Acknowledgments} The authors would like to respectfully express warm thanks to prof. dr. Constantin Udriste, for the inspired suggestion to approach the mechanical problem via optimal control techniques and also for his guidance and valuable remarks, which led to the improvement of this paper.


\end{document}